\documentclass[12pt]{amsart}%
\usepackage{amsmath}
\usepackage{graphicx}
\usepackage{amscd}
\usepackage{geometry}
\usepackage{amsfonts}
\usepackage{amssymb}
\usepackage[mathscr]{eucal}%
\newtheorem{Question}{Question}
\theoremstyle{plain}

\numberwithin{equation}{section}

\begin{document}

\title[ ]{Is a monotone union of contractible open sets contractible?}

\author{Fredric D. Ancel}

\address{Department of Mathematical Sciences, University of Wisconsin - Milwaukee, Box 413, Milwaukee, WI 53201-0413}

\email{ancel@uwm.edu}

\author{Robert D. Edwards}

\address{Department of Mathematics, UCLA, Box 951555 Los Angeles, CA 90095-1555}

\email{rde@math.ucla.edu}

\thanks{}

\date{\today}

\begin{abstract}
This paper presents some partial answers to the following question.
\\[6pt]
\textbf{Question.} If a normal space $X$ is the union of an increasing sequence of open sets $U_1 \subset U_2 \subset U_3 \subset \dots$ such that each $U_n$ contracts to a point in $X$, must $X$ be contractible?
\\[6pt]
The main results of the paper are:
\\[6pt]
\textbf{Theorem 1.} If a normal space $X$ is the union of a sequence of open subsets $\{U_n\}$ such that $cl(U_n) \subset U_{n+1}$ and $U_n$ contracts to a point in $U_{n+1}$ for each $n \geq 1$, then $X$ is contractible.
\\[6pt]
\textbf{Corollary 2.}  If a locally compact $\sigma$-compact normal space $X$ is the union of an increasing sequence of open sets $U_1 \subset U_2 \subset U_3 \subset \dots$ such that each $U_n$ contracts to a point in $X$, then $X$ is contractible.
\end{abstract}
\maketitle

\section{Introduction}

\renewcommand{\thefootnote}
{\ensuremath{\fnsymbol{footnote}}}
\setcounter{footnote}{2}

\indent In 1935, J. H. C. Whitehead, to illustrate a flaw in his own proposed proof of the Poincar\'{e} Conjecture, constructed a contractible open\footnote{A manifold is called $open$ if it is non-compact and has empty boundary.} 3-manifold without boundary that is not homeomorphic to $\mathbb{R}^3$ \cite{7}.  Subsequently it was shown that in each dimension $n \geq 3$, there exist uncountably many non-homeomorphic contractible open $n$-manifolds.  (See \cite{5}, \cite{1} and \cite{3}.)  These spaces illustrate the richness of the topology of manifolds in dimensions greater than 2. 

\indent Proofs that a construction yields a contractible open $n$-manifold that is not homeomorphic to $\mathbb{R}^n$ characteristically have two steps.  First they establish that the constructed space is contractible.  Second they show that it is not homeomorphic to $\mathbb{R}^n$.  While the second step is usually the more interesting and delicate of the two, in this article we focus on methods used to take the first step. 

\indent Typical constructions of contractible open manifolds produce a space $X$ that is the union of an increasing sequence of open subsets $U_1 \subset U_2 \subset U_3 \subset \dots$ such that each $U_n$ contracts to a point in $X$.  With this information one can justify the contractibility of $X$ in various ways.  For instance, if $X$ is a CW complex, then one can observe that all the homotopy groups of $X$ vanish and use a theorem of J. H. C. Whitehead (Corollary 24 on page 405 of \cite{6}) to conclude that $X$ is contractible.  If a more elementary justification is sought which avoids assuming that the space $X$ is a CW complex or appealing to the theorem of Whitehead, then the following theorem provides an approach. 
\\\\\
\noindent \textbf{Theorem 3.} If a normal space $X$ is the union of a sequence of open subsets $\{U_n\}$ and there is a point $p_0 \in U_1$  such that for each $n \geq 1, cl(U_n) \subset U_{n+1}$ and $U_n$ contracts to $p_0$ in $U_{n+1}$ fixing $p_0$, then $X$ is contractible. 
\\\\\
\indent The proof of Theorem 3 is elementary and well known. Observe that Theorem 3 follows immediately from Theorem 1. (Also the first half of the proof of Theorem 1 given below is essentially a proof of Theorem 3. A parenthetical comment in the proof of Theorem 1 marks the point at which the proof of Theorem 3 is complete.) Applying Theorem 3 directly to a space $X$ requires some care in the construction of $X$ to insure that each $U_n$ contracts to an initially specified point $p_0$ in $U_{n+1}$ fixing the point $p_0$.  The motivation behind this paper is to show that we can weaken the hypotheses of Theorem 3 to those of Theorem 1 and thereby remove the requirement that the homotopy contracting $U_n$ to a point in $U_{n+1}$ fixes any particular point.  As a consequence, in the construction of a contractible open manifold, the argument that the constructed object is contractible becomes easier while still relying on principles that are valid in a very broad setting (the realm of normal spaces). 

\indent We remark that the hypothesis that the homotopy contracting $U_n$ to a point in $U_{n+1}$ fix the point can't be dropped with impunity because there exist contractible metric spaces that can't be contracted to a point fixing that point.  The \textit{line of Cantor fans} is a simple non-compact example of such a space.  This space is the countable union $\cup_{n \in \mathbb{Z}}K _n$   in which $K_n$ is the cone in the plane with vertex $(n,0)$ and base $\{n+1\} \times C$ where $C$ is the standard middle-thirds Cantor set in $[0,1]$.  A more complex compact example is the \textit{Cantor sting ray} described in \cite{2}.  (A comparable complete description of the Cantor sting ray can be found in Exercise 7 on pages 18-19 in \cite{4}.)

\indent Although the requirement that the contracting homotopies fix a point can't be omitted without consequence, it is known that it can be omitted if one is willing to impose additional conditions on $X$ as in the following result.
\\\\
\textbf{Theorem 4.} If a normal space $X$ is the union of a sequence of open subsets 
$\{U_n\}$ such that for each $n \geq 1$, $cl(U_n) \subset U_{n+1}$ and $U_n$ contracts to point in $U_{n+1}$, then $X$ is contractible provided that it satisfies the following additional condition.\\[6pt]
$( \ast )$ There is an open subset $V$ of $X$ that contracts to a point $p_0 \in V$ in $X$ fixing $p_0$.
\\\\
\indent Theorem 4 follows from Theorem 3 and the following lemma.
\\\\
\textbf{Lemma 5.} If $W \subset U_1 \subset U_2$ are open subsets of a completely regular space $X$ and if $W$ contracts to a point $p_0 \in W$ in $U_1$ fixing $p_0$ and $U_1$ contracts to a point in $U_2$, then $U_1$ contracts to $p_0$ in $U_2$ fixing $p_0$.
\\\\
\indent Although the proof of Lemma 5 is known and is similar to the proofs of Theorem 1.4.11 on pages 31 and 32 and Exercise 1.D.4 on page 57 of \cite{6}, we follow the referee's recommendation that we include a proof.  

\begin{proof}[Proof of Lemma 5]
There are homotopies $f : W \times [0,1] \rightarrow U_1$ and $g : U_1 \times [0,1] \rightarrow U_2$ such that $f$ contracts $W$ to $p_0$ in $U_1$ fixing $p_0$ and $g$ contracts $U_1$ to a point $q_0$ in $U_2$.
\\\\
\indent \textbf{Step 1.} There is a map $\phi : W \times ([0,1]^2) \rightarrow U_2$ with the following properties.  For all $(x,(s,t)) \in W \times ([0,1]^2)$: $\phi(x,(s,0)) = g(x,s)$, $\phi(x,(0,t)) = x$, $\phi(x,(1,t)) = q_0$ for $0 \leq t \leq 1/2$, $\phi(x,(1,t)) = g(p_0,2-2t)$ for $1/2 \leq t \leq 1$ and $\phi(x,(s,1)) = f(x,s)$.  Observe that $\phi(p_0,(s,1)) = p_0$ for $0 \leq s \leq 1$.

\indent To construct $\phi$ invoke the Tietze Extension Theorem to obtain maps $\lambda, \mu : [0,1]^2 \rightarrow [0,1]$ satisfying the following conditions: $\lambda$ maps $([0,1] \times \{0\}) \cup (\{0\} \times [0,1])$ to $0$, $\lambda$ maps $(\{1\} \times [1/2,1])$ to $1$, $\lambda(s,1) = s$ for $0 \leq s \leq 1$, $\lambda(1,t) = 2t$ for $0 \leq t \leq 1/2$, $\mu$ maps $(\{0\} \times [0,1]) \cup ([0,1] \times \{1\})$ to $0$, $\mu$ maps $(\{1\} \times [0,1/2])$ to $1$, $\mu(s,0) = s$ for $0 \leq s \leq 1$ and $\mu(1,t) = 2-2t$ for $1/2 \leq t \leq 1$. Then define $\phi(x,(s,t)) = g(f(x,\lambda(s,t)),\mu(s,t))$ for all $(x,(s,t)) \in W \times ([0,1]^2)$.
\\\\
\indent \textbf{Step 2.} Let $B = ([0,1] \times \{0\}) \cup (\{0,1\} \times [0,1])$. Observe that $\phi$ can be extended to a map $\psi : (W \times ([0,1]^2)) \cup (U_1 \times B) \rightarrow U_2$ such that for all $x \in U_1$: $\psi(x,(s,0)) = g(x,s)$ for $0 \leq s \leq 1$, $\psi(x,(0,t)) = x$ for $0 \leq t \leq 1$, $\psi(x,(1,t)) = q_0$ for $0 \leq t \leq 1/2$ and $\psi(x,(1,t)) = g(p_0,2-2t)$ for $1/2 \leq t \leq 1$. Observe that $\psi(x,(1,1)) = p_0$ for all $x \in U_1$ and $\psi(p_0,(s,1)) = p_0$ for $0 \leq s \leq 1$.
\\\\
\indent \textbf{Step 3.}  There is a map $r : U_1 \times [0,1]^2 \rightarrow (W \times [0,1]^2) \cup (U_1 \times B)$ that restricts to the identity on $(U_1 \times B) \cup (\{p_0\} \times [0,1]^2)$.  To construct $r$, we exploit the fact that $B$ is a strong deformation retract of $[0,1]^2$.  Hence, there is a homotopy $\delta : [0,1]^2 \times [0,1] \rightarrow [0,1]^2$ that joins the identity on $[0,1]^2$ to a retraction of $[0,1]^2$ onto $B$ while keeping the points of $B$ stationary. Also we exploit the fact that $X$ is a completely regular space to obtain a map $\nu : X \rightarrow [0,1]$ such that $\nu(p_0) = 0$ and $\nu(X-W) = \{1\}$.  Now we define the map $r$ by $r(x,(s,t)) = (x,\delta((s,t),\nu(x)))$ for $(x,(s,t)) \in U_1 \times [0,1]^2$.
\\\\
\indent \textbf{Step 4.} Finally we define the homotopy $\omega : U_1 \times [0,1] \rightarrow U_2$ by $\omega(x,s) = \psi\circ r(x,(s,1))$ for $(x,s) \in U_1 \times [0,1]$.  Then $\omega$ contracts $U_1$ to $p_0$ in $U_2$ fixing $p_0$.
\end{proof}

\begin{proof}[Proof of Theorem 4]
To prove Theorem 4 from Theorem 3 and Lemma 5, observe that under the hypotheses of Theorem 4, there is an $m \geq 1$ for which $p_0 \in U_m$.  Then a neighborhood $W$ of $p_0$ in $V$ can be chosen so that the homotopy contracting $V$ to $p_0$ in X fixing $p_0$ restricts to a homotopy contracting $W$ to $p_0$ in $U_m$.  Then Lemma 5 implies that for each $n \geq m$, $U_n$ contracts to $p_0$ in $U_{n+1}$ fixing $p_0$.  We can now invoke Theorem 3 to conclude that $X$ is contractible.
\end{proof}

\indent Since every manifold and, more generally, every absolute neighborhood retract 
satisfies hypothesis $( \ast )$ of Theorem 4, we have the following corollary.
\\\\
\textbf{Corollary 6.} If an absolute neighborhood retract $X$ is the union of a sequence of open subsets $\{U_n\}$ such that for each $n \geq 1$, $cl(U_n) \subset U_{n+1}$ and $U_n$ contracts to point in $U_{n+1}$, then $X$ is contractible.
\\\\
\indent Observe that Theorem 1 reaches the same conclusion as Theorem 4 and Corollary 6 without assuming hypotheses like $( \ast )$ or that $X$ is an absolute neighborhood retract.  Establishing that the contractibility of $X$ can be proved without imposing such additional restrictions on $X$ is one of the objectives of this paper. 

\indent The authors with to express their appreciation to the Workshop in Geometric Topology for providing a venue and a table full of willing participants - faculty and students - to bat around the questions that gave rise to this article.\\

\section{Proofs of Theorem 1 and Corollary 2}

\begin{proof}[Proof of Theorem 1]
By hypothesis, for each $n \geq 1$, there is a homotopy $f_n : U_{3n} \times [0,1] \rightarrow U_{3n+1}$ such that $f_n(x,0) = x$ for each $x \in U_{3n}$ and $f_n(U_{3n} \times \{1\}) = \{p_n\}$ for some point $p_n \in U_{3n+1}$. We modify each $f_n$ to get a homotopy with domain $X \times [0,1]$ by invoking Urysohn's Lemma to obtain a map $\lambda_n : X \rightarrow [0,1]$ such that $\lambda_n(cl(U_{3n-2})) = \{1\}$ and $\lambda_n(X - U_{3n-1}) = \{0\}$.  
Then we define the homotopy $g_n : X \times [0,1] \rightarrow X$ by:
\[ 
	g_n(x,t) =
	\begin{cases}
		f_n(x,\lambda_n(x)t)             &\text{for $(x,t) \in U_{3n} \times [0,1]$}\\
		x				           &\text{for $(x,t) \in (X - U_{3n-1}) \times [0,1]$}
	\end{cases}
\]
\noindent Therefore:
\begin{itemize}
\item $g_n(x,0) = x$ for each $x \in X$,
\item $g_n(U_{3n-2} \times \{1\}) = \{p_n\}$,
\item $g_n(x,t) = x$ for each $(x,t) \in (X - U_{3n-1}) \times [0,1]$, and
\item $g_n(U_{3n+1} \times [0,1]) \subset U_{3n+1}$.
\end{itemize}

\indent Next we define a map $h : X \times [0,\infty) \rightarrow X$ by stacking the $g_n$'s.  First we define $\phi_n : X \rightarrow X$ by $\phi_n(x) = g_n(x,1)$ for each $x \in X$.  Then we define $h : X \times [0,\infty) \rightarrow X$ by:
\[
	h(x,t) =
	\begin{cases}
		g_1(x,t) 		&\text{for $(x,t) \in X \times [0,1]$}\\
		g_{n+1}(\phi_n\circ\dots\circ\phi_1(x),t - n)  &\text{for $(x,t) \in X \times [n,n+1]$ and $n \geq 1$}
	\end{cases}
\]
\indent Observe that $h(x,0) = g_1(x,0) = x$ for each $x \in X$.  Let 
\begin{align*}
A  = \bigcup_{n=1}^{\infty}(U_{3n-2} \times [n,n+1]).
\end{align*}

For $n \geq 1$, since $\phi_i(U_{3n-2}) \subset U_{3n-2}$ for $1 \leq i \leq n$ and $\phi_n(U_{3n-2}) = \{p_n\}$, it follows that $h(U_{3n-2} \times \{t\}) = \{g_{n+1}(p_n,t - n)\}$ for $t \in [n,n+1]$.  Thus, $h$ is constant on each horizontal slice of $A$; in other words, $h(A \cap (X \times \{t\}))$ is a one-point set for each $t \in [1,\infty)$.  

\indent In this situation, if it were the case that the map $t \mapsto h(A \cup (X \times \{t\})) : [1,\infty) \rightarrow X$ converges to a point $p$ of X as $t$ approaches $\infty$, it follows that one could extend $h$ to a map of 
$X \times [0,\infty]$ to $X$ which would contract $X$ to $p$.  (If we were proving Theorem 3, then this map could be contrived to be the constant map with value $p_0$, and the proof of Theorem 3 would be finished at this point.)  However, in the current situation, there is no reason to expect such 
convergence.  Instead, we introduce another device to establish the contractibility of $X$.  
The exposition of this device is the main contribution of this article.

\indent Let $\sigma : X \rightarrow [3,\infty)$ be a map whose graph is contained in $A$. (See Figure 1.)  One scheme for constructing $\sigma$ is the following.  For each $n \geq 1$, invoke Urysohn's Lemma to obtain a map $\sigma_n : X \rightarrow [0,1]$ so that $\sigma_n(cl(U_{3n-2})) = \{0\}$ and $\sigma_n(X - U_{3n+1}) = \{1\}$.  Then define $\sigma$ by the formula $\sigma(x) = 3 + \sum_{n=1}^{\infty}\sigma_n(x)$.

\begin{figure}[ht]
\centerline{
\includegraphics{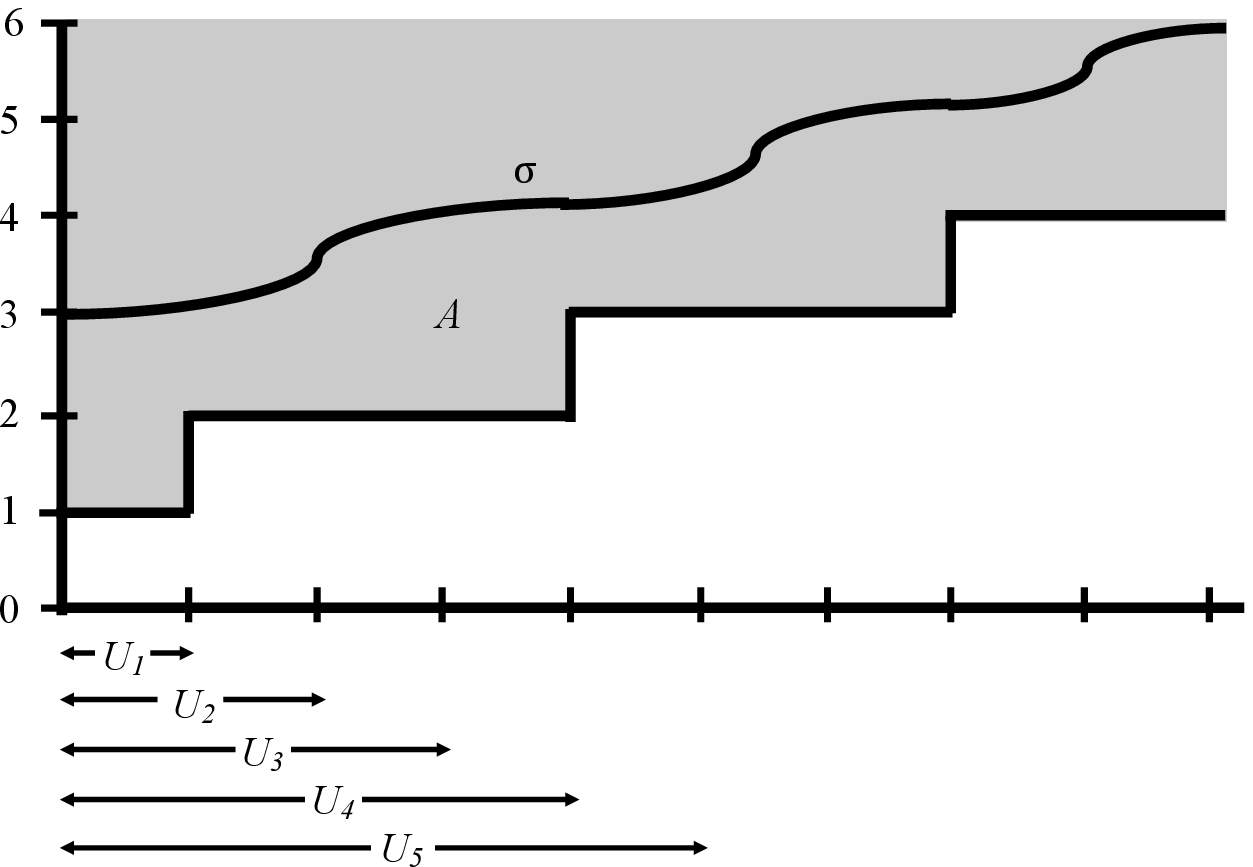}
}
\caption{}
\end{figure}

\indent Define the map $\psi : X \rightarrow X$ by $\psi(x) = h(x,\sigma(x))$. Then a homotopy $k : X \times [0,1] \rightarrow X$ that joins $id_X$ to $\psi$ is defined by $k(x,t) = (x,\sigma(x)t)$. The virtue of $\psi$ that makes it useful at this juncture is that it factors through $[3,\infty)$. To verify this claim, we choose a point $q \in U_1$ and define the map $\tau : [3,\infty) \rightarrow X$ by $\tau(t) = h(q,t)$. Observe that since $\{q\} \times [3,\infty) \subset A$ and $h$ is constant on each horizontal slice of $A$, it follows that $\psi(x) = h(x,\sigma(x)) = h(q,\sigma(x)) = \tau(\sigma(x))$. Thus, $\psi = \tau\circ\sigma$. Therefore, $id_X$ is homotopic to $\psi$ and $\psi$ is null-homotopic because it factors through the contractible space $[3,\infty)$.  Hence, $id_X$ is null-homotopic.  Consequently, $X$ is contractible.
\end{proof}

\begin{proof}[Proof of Corollary 2] Since $X$ is locally compact, every compact subset of $X$ is contained in an open subset of $X$ with compact closure.  Since $X$ is $\sigma$-compact, it follows that $X$ is the union of a sequence of open subsets $\{V_m\}$ such that each $V_m$ has compact closure.  

\indent We will now construct by induction a sequence $\{W_k\}$ of open subsets of $X$ such that for each $k \geq 1$, $V_k \subset W_k$ and $cl(W_k)$ is compact and contracts to a point in $W_{k+1}$.  Begin by letting $W_1 = V_1$.  Next let $k \geq 1$ and assume $W_k$ is an open subset of $X$ such that $V_k \subset W_k$ and $cl(W_k)$ is compact.  Since $\{U_n\}$ is an increasing open cover of the compact set $cl(W_k)$, it follows that $cl(W_k)$ is a subset of some $U_n$.  Therefore, there is a homotopy which contracts $cl(W_k)$ to a point in $X$.  Since the image of this homotopy is compact, there is an $m \geq k+1$ such that the image of this homotopy lies in $V_1 \cup V_2 \cup \dots \cup V_m$.  Let $W_{k+1} = V_1 \cup V_2 \cup \dots \cup V_m$.  Then $V_{k+1} \subset W_{k+1}$, $cl(W_{k+1})$ is compact and $cl(W_k)$ contracts to a point in $W_{k+1}$.  This completes the construction of the $W_k$.

\indent Since $\{V_k\}$ covers $X$ and $V_k \subset W_k$ for each $k \geq 1$, $\{W_k\}$ covers $X$.  Since $cl(W_k)$ contracts to a point in $W_{k+1}$ for each $k \geq 1$, it follows that the hypotheses of Theorem 1 are satisfied.  Hence, $X$ is contractible. \end{proof}

\section{The Question}

\indent We are unsure whether the restriction to locally compact $\sigma$-compact spaces in the hypothesis of Corollary 2 is necessary.  Moreover, there are large and important classes of contractible spaces such as infinite dimensional Hilbert spaces that are neither locally compact nor $\sigma$-compact.  Thus, omitting the hypotheses of local compactness and $\sigma$-compactness from Corollary 2 would greatly broaden its applicability. This brings us back to our original question.

\begin{Question} If a normal space $X$ is the union of an increasing sequence of open sets $U_1 \subset U_2 \subset U_3 \subset \dots$ such that each $U_n$ contracts to a point in $X$, must $X$ be contractible? \end{Question}

\end{document}